\documentclass[10pt]{amsart}
\usepackage{amssymb,latexsym}

\newtheorem{theorem}{Theorem}

\newtheorem{lemma}{Lemma}

\theoremstyle{definition}

\theoremstyle{remark}

\numberwithin{equation}{section}

\newcommand{\sgn}{\text{sgn}}

\newcommand{\R}{\mathbb R}

\begin{document}
\title[ $g$ZK  in weighted Sobolev spaces]{
Well-posedness for the two dimensional generalized Zakharov-Kuznetsov equation in anisotropic  weighted Sobolev spaces}
\author{G. Fonseca}
\address[G. Fonseca]{Departamento  de Matem\'aticas\\
Universidad Nacional de Colombia\\ Bogot\'a\\Colombia}
\email{gefonsecab@unal.edu.co}
\author{M. Pach\'on  }
\address[M. Pach\'on]{Departamento  de Matem\'aticas\\
Universidad Central\\ Bogot\'a\\Colombia}
\email{mpachonh@ucentral.edu.co}
\begin{abstract} We consider the well-posedness of the initial value problem associated to the $k$-generalized Zakharov-Kuznetsov equation in fractional weighted Sobolev spaces $H^s(\R^2)\cap L^2((|x|^{2r_1}+|y|^{2r_2})\,dxdy)$, $s, r_1, r_2 \in \R$. Our method of proof is based on the contraction mapping principle and it mainly relies on the well-posedness results recently obtained for this equation in the Sobolev spaces  $H^s(\R^2)$ and  a new pointwise commutator type formula involving the group induced by the linear part of the equation and the fractional anisotropic weights to be considered.  
\end{abstract}
\maketitle

\section{Introduction}\label{S:1}
Our aim is to study  persistence properties of solutions of  the two dimensional $k$-generalized Zakharov-Kuznetsov equation (gZK) 
in fractional weighted spaces. More precisely we consider  the initial value problem (IVP):
\begin{equation}\label{kgKdV}
\begin{cases}
\partial_tu +\partial_x\Delta u+ u^k\partial_xu=0,\;\;\;\;t\in\R,\,(x,y)\in\R^2,\;\;k\in\mathbb Z^+,\\
u(x,y,0)=u_0(x,y).
\end{cases}
\end{equation}

Let us introduce the weighted Sobolev spaces of our interest
\begin{equation}
\label{spaceZ}
Z_{s,(r_1,r_2)} = H^s(\R^2) \cap L^2((\,|x|^{2r_1}+|y|^{2r_2})\,dxdy),\,\,\, s, r_1, r_2 \in \R.
\end{equation}
We want to show  that for initial data  in this function space the associated IVP is locally well-posed and with some additional assumptions it turns out to be globally well-posed. In general an IVP is said to be locally well-posed (LWP) in a function space $X$ if for each $u_0\in X$ there exist $T>0$ and a unique solution $u\in C([-T,T]:X)\cap\dots=Y_T$ of the equation, such that the map data $\to$ solution is locally continuous from $X$ to $ Y_T$. 

This notion of LWP includes the \lq\lq  persistence" property, i.e. the solution describes a continuous curve on $X$. In particular, this implies that the solution flow of the considered equation defines a dynamical system in $X$.   Whenever $T $ can be taken arbitrarily large we say  that the IVP  is globally well-posed (GWP).

It is important to mention that this family of dispersive equations include the Zakharov-Kuznetsov (ZK) equation ($k=1$) and the modified  Zakharov-Kuznetsov (mZK) equation ($k=2$) which are considered two dimensional versions of the famous Korteweg-de Vries (KdV) and modified Kortewg-de Vries (mKdV) equations respectively. The ZK equation was introduced by Zakharov and Kusnetsov in \cite{ZK}  in the context of plasma physics in order to model the propagation of ion-acoustic waves in magnetized plasma,  for a rigorous proof of this fact see \cite{LLS}. On the other hand mKdV equation is used to describe the propagation of Alfv\'en waves at a critical angle to an undisturbed magnetic field (see \cite{KaO}) and mZK is related to the same type of phenomena  in two dimensions (see \cite{SiBe}). 

In order to motivate our results we remark that for the gKdV IVP:
\begin{equation}\label{gKdV}
\begin{cases}
\partial_tu +\partial_x^3 u+ u^k\partial_xu=0,\;\;\;\;t,\,x\in\R,\;\;k\in\mathbb Z^+,\\
u(x,0)=u_0(x),
\end{cases}
\end{equation}

Kato   showed in \cite{Ka1} the persistence of solutions in the weighted Sobolev spaces 

$$
Z_{s,m}=H^s(\R)\cap L^2(\,|x|^{2m}dx),\;\;\;\;\;s\geq 2m,\;\;\;\;\;m=1,2,\dots
$$
The proof of this result  is based on the commutative property of the operators
\begin{equation}
\label{op-ga}
\Gamma =x-3t\partial_x^2,\;\;\;\;\;\;\;\mathcal L=\partial_t+\partial_x^3,\;\;\;\;\;\text{i.e.}\;\;\;\;[\Gamma;\mathcal L]=0.
\end{equation}
Let us consider the linear IVP 
 \begin{equation}
\label{lKdV}
\begin{cases}
\partial_tv +\partial_x^3 v=0,\;\;\;\;t,\,x\in\R,\\
v(x,0)=v_0(x),
\end{cases}
\end{equation}
and let us denote by $\,\{U(t) \,:\,t\in\R\}$  the unitary group of operators describing its solution, that is
\begin{equation}\label{de-gr}
U(t)v_0(x)= (e^{it\xi^3} \widehat v_0)^{\lor}(x)
\end{equation}
 Then  \eqref{op-ga} is equivalent to 
 \begin{equation}
 \label{formula1}
 x \,U(t) v_0(x) = U(t)(xv_0)(x) +3t U(t) (\partial_x^2 v_0)(x).
 \end{equation}
 
This equality clearly suggests that regularity and decay are strongly related and furthermore in order to obtain persistent properties for the flow in \eqref{lKdV} in those weighted spaces $Z_{s,r}$,  at least twice of  the decay rate $r$ is expected  to be required in regularity, that is, $s\geq 2r$. 

Notice that Kato's result strongly indicates us that this condition should hold even for the non-linear associated IVP. In fact, this was recently proved by Fonseca, Linares and Ponce  in \cite{FLPg} where they extended 
\eqref{formula1} to fractional powers of $|x|$ with the help of a point-wise version of the homogeneous derivative of order $s$ introduced by Stein \cite{St1}.    In that way, they improved Kato's  results for the gKdV in those fractional Sobolev weighted spaces. Their argument, via contraction principle, also required some detail on previous results on LWP and GWP results for the gKdV IVP on the classical Sobolev spaces $H^s(\R)$ obtained by Kenig, Ponce and Vega, see \cite{KPV}, \cite{KPV1}.
     
In view of the ideas just presented in the case of the gKdV equation, the first piece in our analysis is related to some of  the existent theory on LWP and GWP for the gZK equation  on classical Sobolev spaces $H^s(\R^2),$ see remark b) below.
Let us define the regularity index $s_k$ by: 
\begin{equation}
\label{reg}
s_k=
\begin{cases}
&\;3/4\;\,\;\;\;\;\; \text{if} \;\;\;1\leq k\leq 7,\\
&\;1-2/k \;\;\;\;\text{if}\;\;\;\;\;k\geq 8.
\end{cases}
\end{equation}
We now state the results by Linares and Pastor \cite{LP}, \cite{LP2} and Farah, Linares, and Pastor \cite{FaLP} and notice that  some detail of their proof will be included in Section \ref{S:2}.  

\begin{theorem}$($\cite{LP}, \cite{LP2}, \cite{FaLP}$)$\label{th1}. For any $u_0\in H^{s}(\R^2), s>s_k$, there exist $T=T(\|u_0\|_{H^s})>0,$
an space $X_T \subset C([0,T]: H^{s}(\R^2))$ and a unique solution $u\in X_T$ of the IVP \eqref{kgKdV} defined in $[0,T]$. Moreover for any $T'\in(0,T)$ there exists a neighborhood $V$ of $u_0$ in $ H^{s}(\R^2)$ such that the map $\widetilde u_0\to\widetilde u(t)$ from $V$ into $X_{T'}$ is smooth.
\end{theorem}

\underline{Remarks}: (a) The estimate for the length of the time interval of existence with respect of the size of the initial data in $H^s(\R^2)$ can be explicitly obtained  in the proof of Theorem \ref{th1}.

(b) The critical index for the gZK equation \eqref{kgKdV}  turns out to be $s_{c,k}=1-\frac2k$ which can easily be computed by an scaling argument, therefore  it coincides with the regularity index $s_k$ in \eqref{reg}  within the range $ k\geq 8$. Hence we have that for $ k\geq 8$ these results are  optimal . Actually, in \cite{FaLP} it was proven that the  gZK IVP is ill-posed for $s=s_{c,k}$ in the sense that the map data to solution is not uniformly continuous so other approach different from contraction arguments is required in order to lower the LWP regularity. However, for the ZK equation, $k=1$, Gr\"unrock and Herr  in \cite{GH} and Molinet and Pilod in \cite{MP} were able to show LWP in $H^s(\R^2)$ for $s>1/2$  in the context of Bourgain spaces $X^{s,b}$, see \cite{Bo1}. 
Also, recently Ribaud and Vento in \cite{RV} showed local well-posedness  in $H^s(\R^2)$ for $s>1/4$ if $k=2$, $s>5/12$ if $k=3$ and $s>s_{c,k}$ if $k\geq 4$ by working in some Besov spaces. Notice that for  $ 1\leq k\leq 3$ there is still a gap to be filled in the expected LWP theory.  

(c) Regarding GWP results, it is important to mention that for the ZK equation local solutions can be globally defined in $H^1$ with the help of the conserved energy and a Gagliardo-Nirenberg inequality. In the case of the gZK equation, global solutions in $H^1$, and even in a larger space in the case of the mZK, $k=2$, are obtained if in addition it is assumed that the initial data is small enough, see \cite{LP}, \cite{LP2} and \cite{FaLP}.
\vskip .2in		
 Next,  let us explicitly introduce the group associated to the linear ZK equation:

\begin{equation}\label{grZK}
W(t)v_0(x,y)= (e^{it(\xi^3+\xi \eta^2)} \widehat v_0)^{\lor}(x,y).
\end{equation}

Following the strategy used to deal  with the gKdV equation in weighted spaces in \cite{FLPg}, our second task is directed to get an extension of formula \eqref{formula1} but with the linear group in \eqref{grZK} instead. More precisely we have our first result:
 
\vskip.1in
 \begin{theorem} \label{point-wise}
Let $r_1, r_2\in (0,1)\,\, , s\geq 2\max\{r_1, r_2\}$ and $\{W(t)\,:\,t\in\R\}$ be the unitary group of operators defined in \eqref{grZK}.   If
\begin{equation} \label{hyp1}
u_0\in  Z_{s,r}=H^{}(\R^2)\cap L^2((\,|x|^{2r_1}+|y|^{2r_2})\,dxdy),
\end{equation}
then for all $t\in \R$ and for almost every $(x,y)\in \R^2$
\begin{equation}
\label{comx}
|x|^{r_1} W(t) u_0(x,y) = W(t) (|x|^{r_1} u_0)(x,y) + W(t)\{\Phi_{1,t,r_1}(\widehat u_0)(\xi,\eta)\}^{\lor}(x,y)
\end{equation}
 with
 \begin{equation}
\label{cmx-norm}
\| \{ \Phi_{1, t,r_1}(\widehat u_0)(\xi,\eta)\}^{\lor}\|_2\leq c(1+|t|)(\|u_0\|_2 + \| D_x^{s}u_0\|_2  +  \| D_y^{s}u_0\|_2)
\end{equation}
and
\begin{equation}
\label{cmy}
|y|^{r_2} W(t) u_0(x,y) = W(t) (|y|^{r_2} u_0)(x,y) + W(t)\{\Phi_{2,t,r_2}(\widehat u_0)(\xi,\eta)\}^{\lor}(x,y)
\end{equation}
 with
 \begin{equation}
\label{cmy-norm}
\| \{ \Phi_{2,t,r_2}(\widehat u_0)(\xi,\eta)\}^{\lor}\|_2\leq c(1+|t|)(\|u_0\|_2 + \| D_x^{s}u_0\|_2  +  \| D_y^{s}u_0\|_2).
\end{equation}

Moreover, if in addition to \eqref{hyp1} we suppose  that for $\beta\in (0,\min\{r_1, r_2\})$
\begin{equation}
\label{hyp2}
D^{\beta} (|x|^{r_1}u_0), D^{\beta} (|y|^{r_2}u_0)\in  L^{2}(\R^2)\;\;\;\;\;\text{and}\;\;\;\;\;u_0\in H^{\beta+s}(\R^2),
\end{equation}
 then for all $t\in \R$ and for almost every $(x,y)\in \R^2$
\begin{equation}
\label{d21}
\begin{aligned}
&D^{\beta} (|x|^{r_1} W(t) u_0)(x,y) \\
\\
&= W(t) (D^{\beta}|x|^{r_1} u_0)(x,y) + W(t)(D^{\beta} (\{ \Phi_{1,t,r_1}(\widehat u_0)(\xi,\eta)\}^{\lor}))(x,y)
\end{aligned}
\end{equation}
and
\begin{equation}
\label{d22}
\begin{aligned}
&D^{\beta} (|y|^{r_2} W(t) u_0)(x,y) \\
\\
&= W(t) (D^{\beta}|y|^{r_2} u_0)(x,y) + W(t)(D^{\beta} (\{ \Phi_{2,t,r_2}(\widehat u_0)(\xi,\eta)\}^{\lor}))(x,y)
\end{aligned}
\end{equation}
  with
 \begin{equation}
\label{d212-norm}
\|D^{\beta}(\{ \Phi_{j,t,r_2}(\widehat u_0)(\xi,\eta)\}^{\lor})\|_2\leq c(1+|t|)(\|u_0\|_2 + \| D_x^{\beta+s}u_0\|_2 + \| D_y^{\beta+s}u_0\|_2), 
\end{equation} 
for $j=1,2.$
\end{theorem}

\underline{Remarks}: (a) As we mentioned above, this type of formula was recently established by Fonseca, Linares and Ponce \cite{FLPg} in the context of the Airy group and more generally it also holds for the group associated to the dispersion generalized Benjamin-Ono equation:
\begin{equation}\label{DGBO}
\begin{cases}
\partial_t u -   D^{1+a}_x\partial_x u = 0, \qquad t, x\in \R,\;\;\;\; 0\leq a < 1,\\
u(x,0) = u_0(x),
\end{cases}
\end{equation}
where $D^s$ denotes  the homogeneous derivative of order $s\in\R$,
$$
D^s=(-\partial_x^2)^{s/2}\;\;\;\text{so}\;\;\;D^s f=c_s\big(|\xi|^s\widehat{f}\,\big)^{\vee}, \;\;\;\text{with} \;\;\;D^s=(\mathcal H\,\partial_x)^s,
$$
and $\mathcal H$ denotes the Hilbert transform,
\begin{equation*}
\mathcal H f(x)=\frac{1}{\pi}\lim_{\epsilon\downarrow 0}\int\limits_{|y|\ge \epsilon} \frac{f(x-y)}{y}\,dy=(-i\,\sgn(\xi) \widehat{f}(\xi))^{\vee}(x).
\end{equation*}

For the problem we are dealing with we adapt those one dimensional situations and carefully handle estimates in the two Fourier space variables.

(b) The  proof of Theorem \ref{point-wise}  is based on a characterization of the generalized Sobolev space
 \begin{equation}
 \label{def0}
 L^{\alpha,p}(\R^n)= (1-\Delta)^{-\alpha/2} L^{p}(\R^n),\;\;\;\;\alpha\in (0,2),\;\;p\in(1,\infty),
 \end{equation}
 due to E. M. Stein \cite{St1} (see Theorem \ref{th5} below).

Now we state our second result concerning local well-posedness of the gZK equation in weighted spaces:

\begin{theorem}
\label{th3}
Let  $u\in C([0,T]: H^{s}(\R^2)), s>s_k$  denote the local solution of the IVP  \eqref{kgKdV}  provided by Theorem \ref{th1}. Let us assume that 
$(|x|^{r_1}+|y|^{r_2})u_0\in L^2(\R^2)$ with $s$ satisfying $ 0<2\max\{r_1, r_2\} \leq s$, then 
\begin{equation}
\label{class2}
\;\;u\in C([0,T]: Z_{s,(r_1, r_2)}).
\end{equation}

For any $T'\in(0,T)$ there exists a neighborhood $V$ of $u_0$ in
$ H^{s}(\R^2)\cap L^2((|x|^{2r_1}+|y|^{2r_2})dxdy)$ such that the map $\widetilde u_0\to\widetilde
u(t)$ from $V$ into the class defined by $X_T$ in Theorem 1 and \eqref{class2} with $T'$
instead of $T$ is smooth.
\end{theorem}
\underline{Remarks}:  (a) We observe that Theorem \ref{th3} guarantees that the  persistent property in the weighted space 
$Z_{s,(r_1, r_2)}$ holds in the same time interval $[0,T]$  given by Theorem \ref{th1}, where
$T$ depends only on $\,\|u_0\|_{H^s}$. 

(b) It is expected  that the condition $s\geq 2\max\{r_1, r_2\}$ in Theorem \ref{th3} is optimal as it was shown in \cite{ILP} in the case of the gKdV equation.
More precisely, \eqref{class2} would hold only if and only if $s\geq 2\max\{r_1, r_2\}$.

  (c) Notice that for $ k=1$ the LWP results in  \cite{GH} and \cite{MP}  hold in a much larger space involving Bourgain spaces $X^{s,b}, s>1/2$, and similarly for  $2\leq k\leq 7$ LWP results in \cite{RV} involve Besov spaces but so far it is not clear for us how to handle our weights in those spaces. 

(d) It is interesting to mention an important difference in the way we obtain persistent properties in these weighted spaces for dispersive type equations. Theorem \ref{th3} is established via the contraction principle  as it was made for semi-linear  Schr\"odinger, gKdV, regularized Benjamin-Ono and the fifth order KdV equations, see \cite{NP}, \cite{FLPg}, \cite{FRS} and \cite{BJJ} respectively. However,  this technique couldn't be used for the dispersive family  DGBO in \eqref{DGBO}  with the quadratic non-linearity $u\partial_xu$ nor for the famous Benjamin-Ono equation. We point out that even that we continue having at hand Theorem \ref{th1}, the dispersion on these equations is too weak to overcome the nonlinear effects in the  contraction argument. Nevertheless, optimal persistency results in weighted Sobolev spaces were indeed  attained via energy estimates for the associated IVPs, see \cite{FLP} and \cite{FP}. See also \cite{AP} regarding a 2D ZK-BO equation. 

(e) Recently, it was proved in \cite{BJM2} a similar result for isotropic weights for the ZK equation, $k=1$ in (\ref{kgKdV}). Their argument of proof depends upon the range of values for the regularity of the initial data in $H^s(\R^2)$. For those values of $s$ such that $\frac34<s\leq 1$ it relies on a symmetrization argument performed to the ZK equation as in \cite{GH}, the proof of a LWP theory for the resultant symmetrized evolution equation on these weighted Sobolev spaces via a contraction argument and  the help of another characterization of the generalized Sobolev spaces $L^{\alpha,p}(\R^n)$ in (\ref{def0}) obtained  by an alternate  Stein's derivative to the one in Theorem 4 below and used in previous works related for the BO and GDBO equations found in in \cite{FP}, \cite{FLP} respectively. In this case the existence time depends on  the size of the initial data on $Z_{s,r}$. For $s>1$ they use the already LWP theory in \cite{LP}, continuos dependence for the flow and energy type arguments in order to prove persistence for the weights. In this case, the time of existence solely depends on the size of the initial data on $H^s(\R^2)$ as in our results. We consider that our proof for this case is simpler and  that Lemma 1 below which basically establishes a commutator property between weights and groups associated to the linear evolution of dispersive equations can be applied in many other multidimensional models. Last, we could recover the contraction argument to get solutions in a more general setting with anisotropic weights and a wider class of power non-linearities.

(f) As a consequence of Theorem \ref{th1},  its remark (c) and our proof of Theorem \ref{th3}, local results globally extend when $s\geq1$ with arbitrary size of the initial data for k=1 and for small enough initial data for $k\geq3$ ( see \cite{FaLP}). For mZK, $s>53/63$ and appropriate smallness of the initial data guarantee such extension of the local solutions (see \cite{LP2}). 



The paper is organized as follows. In section 2 we introduce Stein's derivatives and some detail on known results on LWP for the gZK equation. The proof of Theorem \ref{point-wise} will be given in Section 3.  In Section 4 we will present the proof of Theorem \ref{th3}.

$\bf{Notations.}$ $\|\cdot\|_{p}$ denotes the norm in the Lebesgue space $L^p(\R^n)$. 

Let $\alpha$ be a complex number, the homogeneous derivatives $D_x^\alpha, D_y^\beta$ for functions in $\R^2$ are defined via Fourier transform by $\widehat{D_x^\alpha f}(\xi,\eta)=|\xi|^\alpha\hat f(\xi,\eta)$ and $\widehat{D_y^\alpha f}(\xi,\eta)=|\eta|^\alpha\hat f(\xi,\eta)$ respectively. 

We consider the Lebesgue space-time $L^r_TL^p_xL_y^q$ spaces with $1\leq p, q, r<\infty$ endowed with the norm

$$
\|f\|_{L^r_TL^p_xL_y^q}=\left(\int_0^T\left(\int^\infty_{-\infty}\left(\int^\infty_{-\infty}|f(x,y,t)|^q dy\right)^{\frac{p}{q}}dx\right)^{\frac{r}{p}}dt\right)^{\frac{1}{r}}
$$   
with the usual modifications when  $p=\infty$ or $q=\infty$ or $r=\infty.$

In general $c$ denotes a universal constant which may change, increase, from line to line.

\section{Preliminary results}\label{S:2}

Let us start  with a characterization of the Sobolev space
 \begin{equation}
 \label{def00}
 L^{\alpha,p}(\R^n)= (1-\Delta)^{-\alpha/2} L^{p}(\R^n),\;\;\;\;\alpha\in (0,2),\;\;p\in(1,\infty),
 \end{equation}
due to E. M. Stein \cite{St1}. For $\alpha \in (0,2)$ define

 \begin{equation}
 \label{def1}
 D_{\alpha}f(x) = \lim_{\epsilon \to 0} \frac{1}{c_{\alpha}} \,\int_{|y| \geq \epsilon}  \frac{f(x+y)-f(x)}{|y|^{n+\alpha}}dy,
 \end{equation}
where $\;c_{\alpha}=\pi^{n/2}\,2^{-\alpha}\,\Gamma (-\alpha/2)/\Gamma((n+2)/2)$.

  As it was remarked in  \cite{St1} for appropriate $f$, for example $f\in \mathcal S(\R^n)$, one has
 \begin{equation}
 \label{pro1}
\widehat{ D_{\alpha}f}(\xi) = \widehat{D^{\alpha} f}(\xi)\equiv  |\xi|^{\alpha} \,\widehat f(\xi).
\end{equation}

 The following  result concerning  the $L^{\alpha,p}(\R^n)= (1-\Delta)^{\alpha/2} L^{p}(\R^n)$ spaces
 was established  in \cite{St1},
\begin{theorem}$($\cite{St1}$)$ \label{th5}
Let $\alpha\in (0,2)$ and $p\in(1,\infty)$. Then $f\in  L^{\alpha,p}(\R^n)$ if and only if
\begin{equation}
\label{d1}
\begin{cases}
&\;(a)\;\, f\in L^p(\R^n),\\
\\
&\;(b)\;\;D_{\alpha} f\in L^p(\R^n),\;\;\;\;\;\;\,\;\;\;\;\;\;\;\;\;\;\;\;(D_{\alpha} f (x)\;\;\text{defined in \eqref{def1}}),
\end{cases}
\end{equation}
with
 \begin{equation}\label{d1-norm}
\|f\|_{\alpha,p}= \|(1-\Delta)^{\alpha/2} f\|_p\simeq \|f\|_p+\|D_{\alpha}    f\|_p\simeq \|f\|_p+\|\,D^{\alpha}       f\|_p.
\end{equation}
 \end{theorem}

 Notice that if   $f, fg: \R^n \to \R \in L^{\alpha,p}(\R^n)$  and $ g\in L^{\infty}(\R^n)\cap C^2(\R^n)$ and consider Stein's derivatives en each $j-$th direction in $\R^n$, Stein's partial derivatives,
we have that

 \begin{equation}
\label{d2}
\begin{aligned}
D_{j,\alpha}(fg)(x) & = \lim_{\epsilon \to 0} \frac{1}{d_{\alpha}} \,\int_{|y| \geq \epsilon}  \frac{f(x+y\, \vec e_j)\,g(x+y\, \vec e_j)-f(x)\,g(x)}{|y|^{1+\alpha}}\,dy\\
\\
&= \lim_{\epsilon \to 0} \frac{1}{d_{\alpha}} \int_{|y| \geq \epsilon}  g(x)\frac{f(x+y\,\vec e_j)-f(x)}{|y|^{1+\alpha}}\,dy\\
\\
&\;\;\;+ \lim_{\epsilon \to 0} \frac{1}{d_{\alpha}} \int_{|y| \geq \epsilon}  \frac{(g(x+y\,\vec e_j)-g(x))f(x+y\,\vec e_j)}{|y|^{1+\alpha}}\,dy\\
\\
&= g(x)\,D_{j,\alpha}f(x) + \Lambda_{j, \alpha}\left((g(\cdot+y\,\vec e_j)-g(\cdot))f(\cdot+y\, \vec e_j)\right)(x).
\end{aligned}
\end{equation}

  In particular, if $ g(x) =e^{i\,t \,\varphi(x)}$, then
  \begin{equation}
\label{d3}
\begin{aligned}
  &\Lambda_{j,\alpha}((g(\cdot+y\,\vec e_j)-g(\cdot))f(\cdot+y\,\vec e_j))(x)\\
  \\
  & =
   \lim_{\epsilon \to 0} \frac{1}{d_{\alpha}} \int_{|y| \geq \epsilon}  \frac{(g(x+y\,\vec e_j)-g(x))f(x+y\,\vec e_j)}{|y|^{1+\alpha}}\,dy\\
   \\
  &= e^{i\,t\,\varphi(x)}\, \lim_{\epsilon \to 0}\, \frac{1}{d_{\alpha}}\int_{|y| \geq \epsilon} \frac{e^{i\,t(\varphi(x+y\,\vec e_j)-\varphi(x))}-1}{|y|^{1+\alpha}}\,f(x+y\,\vec e_j)\, dy \\
  &= e^{i\,t\,\varphi(x)}\, \Phi_{j, \varphi, \alpha}(f)(x).
  \end{aligned}
  \end{equation}

  Thus,  we obtain the identity
  \begin{equation}
\label{d4}
  D_{j,\alpha}(e^{i\,t\,\varphi(\cdot)}\,f)(x)= e^{i\,t\,\varphi(x)}\,D_{j,\alpha} f(x) +e^{i\,t\,\varphi(x)} \,\Phi_{j, \varphi, \alpha}(f)(x).
 \end{equation}

Now we restrict to $\R^2$ and choose as the phase function the one from the group associated to the linear ZK equation in \eqref{grZK}  
 \begin{equation}
\label{d5}
 \varphi(x_1, x_2) = x_1^3 + x_1x_2^2.
\end{equation}
 For $j = 1,2,$ we shall obtain a bound for 
\begin{equation}
\label {goal1}
\|\Phi_{j, \alpha}(f)\|_p =  \left \|\,\lim_{\epsilon \to 0}\int_{|y| \geq \epsilon} \frac{e^{i\,t(\varphi(x+y\,\vec e_j)-\varphi(x))}-1}{|y|^{1+\alpha}}\,f(x+y\,\vec e_j)\, dy \, \right \|_p.
\end{equation}
Indeed this is achieved in our first result
 \begin{lemma}
 \label{Lemma1}
Let $\,\alpha\in (0,1)$,  and $p\in(1,\infty)$. If
$$
f\in  L^{\alpha,p}(\R) \ \ \text{and} \ \ f \in L^p((1 + x_1^2 + x_2^2)^{\alpha p}dx_1\,dx_2),
$$
then for all $t\in \R$ and for almost every $(x_1,x_2)\in \R^2$

\begin{equation}
\begin{aligned}\label{d10}
D_{j,\alpha} (e^{it(x_1^3 + x_1x_2^2)}\,f)(x_1, x_2)
=& e^{it(x_1^3 + x_1x_2^2)}\, D_{j,\alpha}f(x_1, x_2)\\
&+ e^{it(x_1^3 + x_1x_2^2)} \,\Phi_{j,t,\alpha}(f)(x_1, x_2),
\end{aligned}
\end{equation}

 with

\begin{equation}\label{d10a}
\Phi_{1, t, \alpha}(f)(x_1, x_2) = \lim_{\epsilon \to 0}\frac{1}{d_{\alpha}}\int_{|y|\geq \epsilon} \frac{e^{it(\varphi(x_1 + y,x_2)- \varphi(x_1,x_2))}- 1}{|y|^{1 + \alpha}} \, f(x_1 + y, x_2)\, dy
\end{equation}
and
\begin{equation}\label{d10b}
\Phi_{2, t, \alpha}(f)(x_1, x_2) = \lim_{\epsilon \to 0}\frac{1}{d_{\alpha}}\int_{|y|\geq \epsilon} \frac{e^{it(\varphi(x_1,x_2+y)- \varphi(x_1,x_2))}- 1}{|y|^{1 + \alpha}} \, f(x_1, x_2+y)\, dy
\end{equation}
satisfying
\begin{equation}
\label{d10-norm}
\| \Phi_{j,t,\alpha}(f)\|_p\leq c_{\alpha}(1+|t|)( \|f\|_p + \| \,(1 + x_1^2 + x_2^2)^{\alpha} \,f\|_p),
\end{equation}
\end{lemma}
for j=1,2.
\begin{proof}
Since our interest resides in the parameter $\alpha\in (0,1)$, then we are allowed to pass the absolute value inside the integral sign in
 \eqref{goal1}.\\

  At different parts of our work we will make use of either of the elementary estimates
  \begin{equation}
  \begin{aligned}
  \label{d6}
  \begin{cases}
  &(a)\;\;\;\;\forall \,\theta\in \R\;\;\;\;\;\;|e^{i\theta}-1|\leq 2,\\
  &(b) \;\;\;\;\forall \,\theta\in\R\;\;\;\;\;\;|e^{i\theta}-1|\leq |\theta|.
  \end{cases}
  \end{aligned}
  \end{equation}

Let us consider first Stein's derivative with respect to $x_1$, i.e. $j= 1$ in \eqref{d10}. From \eqref{d6} (a) and Minkowski's integral inequality  it follows that
   \begin{equation}
  \label{step1}
  \begin{aligned} 
 &\left \Vert \int_{|y|\geq \frac{1}{100}} \frac{e^{it(\varphi(x_1 + y,x_2)- \varphi(x_1,x_2))}- 1}{|y|^{1 + \alpha}} \, f(x_1 + y, x_2)\, dy \right \Vert_p \\
   &\leq  \int_{|y|\geq \frac{1}{100}} \left \Vert \frac{2|f(x_1 + y, x_2)|}{|y|^{1 + \alpha}} \right \Vert_p \, dy \\
	&\leq c \int_{|y|\geq \frac{1}{100}} \frac{\left \Vert f(x_1 + y, x_2)  \right \Vert_p}{|y|^{1 + \alpha}} \, dy \\
&\leq c\left \Vert f  \right \Vert_p\int_{|y|\geq \frac{1}{100}} \frac{1}{|y|^{1 + \alpha}} \, dy \\
&\leq c_{\alpha}\left \Vert f  \right \Vert_p .
  \end{aligned}
  \end{equation}

  Now let us consider the estimate
   \begin{equation}
\label {goal2}
\left \Vert \int_{\epsilon \leq |y| \leq \frac{1}{100}} \frac{e^{it(\varphi(x_1 + y, \, x_2)- \varphi(x_1, \, x_2))}- 1}{|y|^{1 + \alpha}} \, f(x_1 + y, x_2)\, dy \right \Vert_p.
\end{equation}

Inequality \eqref{d6} (b)  yields
\begin{equation}
  \label{step2}
  \begin{aligned}
   \left |  e^{it(\varphi(x_1 + y, x_2) - \varphi(x_1, x_2))} - 1 \right | &\leq |t(\varphi(x_1 + y, x_2) - \varphi(x_1, x_2))| \\
   &= |t| |y| \left | \int_0^1 \partial_{x_1}\varphi(x_1 + sy,x_2)\, ds \right | . 
 \end{aligned}
\end{equation}

   For $x_1,x_2$ in the ball  $B_{100}(0)=\{(x_1,x_2)/\,x_1^2 + x_2^2 < 100\}$ we obtain

  \begin{equation}
  \begin{aligned}
  	|\partial_{x_1}\varphi(x_1 + sy,x_2)| &= 3(x_1 + sy)^2 + x_2^2  \\
    &\leq 6x_1^2 + 6s^2y^2+x_2^2  \\
    &\leq c,  
  \end{aligned}
  \end{equation}
  
  and therefore \[  \left |  e^{it(\varphi(x_1 + y, x_2) - \varphi(x_1, x_2))}-1 \right | \le c \,|t| \, |y|.  \]

   Hence our estimate is summarized  as:
\begin{equation}
  \label{step3}
  \begin{aligned}
&\left \Vert \int_{\epsilon \leq |y| \leq \frac{1}{100}} \frac{e^{it(\varphi(x_1 + y, \, x_2)- \varphi(x_1, \, x_2))}- 1}{|y|^{1 + \alpha}} \, f(x_1 + y, x_2)\, dy \right \Vert_{L^p(B_{100}(0))} \\
 &\leq  c\int_{|y|\leq \frac{1}{100}} \left \Vert \frac{|t||y||f(x_1 + y, x_2)|}{|y|^{1+\alpha}} \right \Vert_{L^p(B_{100}(0))} \, dy \\
	&\leq c|t|\,\int_{|y|\leq \frac{1}{100}} \frac{1}{|y|^{\alpha}}\left \Vert f(x_1 + y, x_2)  \right \Vert_{L^p(B_{100}(0))} \, dy \\
& \leq c|t|\,\left \Vert f  \right \Vert_p\int_{|y|\leq \frac{1}{100}} \frac{1}{|y|^{\alpha}} \, dy \\
& \leq c_{\alpha}\,|t|\,\left \Vert f  \right \Vert_p. 
 \end{aligned}
 \end{equation}

  From the above estimates we now have to consider in \eqref{goal1} the region:
  \[
  |y|\leq 1/100,\;\;\;\;\;\text{and}\;\;\;\;\;x_1^2 + x_2^2 \geq 100.
  \]

  We sub-divide it into two parts:
 \begin{equation}
  \label{2sets}
  (a)\;\;\;|y| \leq \frac{1}{1 + x_1^2 + x_2^2},\;\;\;\;\;(b)\;\;\;|y| \geq \frac{1}{1 + x_1^2 + x_2^2}.
  \end{equation}  
We first assume $|y| \leq \frac{1}{1 + x_1^2 + x_2^2}$ and  observe that $|\partial_{x_1}\varphi(x_1 + sy, x_2)| \le c(1 + x_1^2 + x_2^2)$ in this region. 
With the help of the  change of variable $\tilde y=(1 + x_1^2 + x_2^2)y$, inequality \eqref{d6} (b), the argument in  \eqref{step2} and  Minkowski's inequality
we obtain 
\begin{align*}
&\left \Vert \int_{ |y|\leq \frac{1}{1 + x_1^2 + x_2^2}}\frac{e^{it(\varphi(x_1 + y, \, x_2)- \varphi(x_1, \, x_2))}- 1}{|y|^{1 + \alpha}} \, f(x_1 + y, x_2)\, dy \right \Vert_{L^{p}(B_{100}(0)^c)} \\
& \leq c_{\alpha} \left \Vert \int_{|y|\leq \frac{1}{1 + x_1^2 + x_2^2}}  \frac{|t||y|(1 + x_1^2 + x_2^2)|f(x_1 + y, x_2)|}{|y|^{1+\alpha}} \, dy \right \Vert_{L^{p}(B_{100}(0)^c)} \\
	&\leq c_{\alpha} \left \Vert \int_{|\tilde{y}|\leq 1}  \frac{|t|(1 + x_1^2 + x_2^2)^{\alpha}|f(x_1 + \frac{\tilde{y}}{1 + x_1^2 + x_2^2}, x_2)|}{|\tilde{y}|^{\alpha}} \, d\tilde{y} \right \Vert_{L^{p}(B_{100}(0)^c)} \\
	&\le c_{\alpha}\left \Vert \int_{|\tilde{y}|\leq 1}  \frac{|t|\left( 1 +  (x_1+ \frac{\tilde{y}}{1 + x_1^2 + x_2^2})^2 + x_2^2\right)^{\alpha}|f(x_1 + \frac{\tilde{y}}{1 + x_1^2 + x_2^2}, x_2)|}{|\tilde{y}|^{\alpha}} \, d\tilde{y} \right \Vert_{L^{p}(B_{100}(0)^c)}  \\
	&+ c_{\alpha}\left \Vert \int_{|\tilde{y}|\leq 1}  \frac{|t|\left(\frac{\tilde{y}}{1 + x_1^2 + x_2^2}\right)^{2\alpha}|f(x_1 + \frac{\tilde{y}}{1 + x_1^2 + x_2^2}, x_2)|}{|\tilde{y}|^{2\alpha}} \, d\tilde{y} \right \Vert_{L^{p}(B_{100}(0)^c)}.	
\end{align*}

Now we perform a second change of variable $u = x_1 + \frac{\tilde{y}}{1 + x_1^2 + x_2^2}, v=x_2 $ and since 

\begin{equation}
\label{changeofvariable}
\frac{\tilde y}{(1+x_1^2 + x_2^2)}=y,\;\;\;\;\; |y|\leq 1/100,\;\;\;\;\;x_1^2 + x_2^2 \geq 100,\;\;\;\; \;\text{so} \;\;\;du \sim dx_1,
 \end{equation}

and we conclude that

\begin{align*}
	&\left \Vert \int_{ |y|\leq \frac{1}{1 + x_1^2 + x_2^2}}\frac{e^{it(\varphi(x_1 + y, \, x_2)- \varphi(x_1, \, x_2))}- 1}{|y|^{1 + \alpha}} \, f(x_1 + y, x_2)\, dy \right \Vert_{L^{p}(B_{100}(0)^c)} \\
		&\leq c_{\alpha} |t| \, \int_{|\tilde{y}|\leq 1}  \left \Vert   \frac{(1 + x_1^2 + x_2^2)^{\alpha}|f(x_1, x_2)|}{|\tilde{y}|^{\alpha}}  \right \Vert_p\, d\tilde{y}   \\
	&+  c_{\alpha} |t| \, \int_{|\tilde{y}|\leq 1}   \left \Vert   \frac{|f(x_1, x_2)|}{|\tilde{y}|^{\alpha}}  \right \Vert_p\, d\tilde{y} \\
	&\le c_{\alpha}\,|t|\,(\Vert f \Vert_p + \Vert(1 + x_1^2 + x_2^2)^{\alpha}f \Vert_p  ).
\end{align*}

Next suppose that $|y| \ge \frac{1}{1 + x_1^2 + x_2^2}$. Changing variable, $\tilde y=(1 + x_1^2 + x_2^2)y$, using \eqref{d6} part (a), Minkowski's inequality, and a second change of
variable as in \eqref{changeofvariable} we get
\begin{align*}
&\left \Vert \int_{\frac{1}{1+x_1^2 + x_2^2}}\frac{e^{it(\varphi(x_1 + y, \, x_2)- \varphi(x_1, \, x_2))}- 1}{|y|^{1 + \alpha}} \, f(x_1 + y, x_2)\, dy \right \Vert_{L^{p}(B_{100}(0)^c)} \\
	&\leq c\left \Vert \int_{\frac{1}{1+x_1^2 + x_2^2} \le |y| \le \frac{1}{100}}  \frac{|f(x_1 + y, x_2)|}{|y|^{1+\alpha}} \, dy \right \Vert_{L^{p}(B_{100}(0)^c)} \\
	&\leq  c\left \Vert \int_{1 \le |\tilde{y}|\le \frac{1 + x_1^2 + x_2^2}{100}}  \frac{(1 + x_1^2 + x_2^2)^{\alpha}|f(x_1 + \frac{\tilde{y}}{1 + x_1^2 + x_2^2}, x_2)|}{|\tilde{y}|^{1+\alpha}} \, d\tilde{y} \right \Vert_{L^{p}(B_{100}(0)^c)} \\
	& \leq c_{\alpha}\left \Vert \int_{|\tilde{y}|\geq 1}  \frac{\left(1 + (x_1+ \frac{\tilde{y}}{1 + x_1^2 + x_2^2})^2 + x_2^2\right)^{\alpha}f(x_1 + \frac{\tilde{y}}{1 + x_1^2 + x_2^2}, x_2)\chi_{A}}{|\tilde{y}|^{1+\alpha}} \, d\tilde{y} \right \Vert_{L^{p}(B_{100}(0)^c)}  \\
	&+ \left \Vert \int_{|\tilde{y}|\geq 1}  \frac{\left(\frac{\tilde{y}}{1 + x_1^2 + x_2^2}\right)^{2\alpha}f(x_1 + \frac{\tilde{y}}{1 + x_1^2 + x_2^2}, x_2)\chi_{A}}{|\tilde{y}|^{1+\alpha}} \, d\tilde{y} \right \Vert_{L^{p}(B_{100}(0)^c)}	\\
	&\le c_{\alpha}(\Vert f \Vert_p + \Vert(1 + x_1^2 + x_2^2)^{\alpha}f \Vert_p ),
\end{align*}
where $A=\{\tilde y : |(x_1,x_2)| \ge \sqrt{100|\tilde y|-1}\}$.
\vskip0.1in
On the other hand, regarding Stein's derivative with respect to $x_2$, i.e. $j=2$ in \eqref{d10}, we notice that the useful inequality, $2x_1x_2\leq x_1^2+x_2^2$, allows us to perform the same computations and obtain exactly the same bound \eqref{d10-norm}. 
\end{proof}
\vskip.1in
Next we focus in the local $H^s$ theory in Theorem \ref{th1} for the gZK IVP in 2D from the works by Linares, Pastor and Farah, see \cite{LP}, \cite{LP2}, \cite{FaLP} and references therein.
Although for every $k=1,2,3, ...$ the proof is accomplished by the contraction principle, the $X_T$ space is different  according to the nonlinearity degree $k$ (different space-time norms are involved in the choice of $X_T$).

Our persistence result in weighted spaces for gZK strongly depends on this theorem so for the sake of clearness we emphasize some of aspects of its proof  for all included nonlinearities  with $k=1,2,3,...$.  

We start by noticing that the method of proof is performed via the Picard iteration applied to the integral version of the gZK IVP given by:

\begin{equation}
\label{in-inhoeqk}
\Psi(u(t))=W(t)u_0-\int_0^{t} W(t-t')(u^{k} \; \partial_x u)(t')dt',\;\;t\in[0,T],
\end{equation} 
where the solution space $X_T \subset C([0,T]: H^{s}(\R^2))$ is determined by the norms
\vskip.1in
\begin{itemize}
\item  For $k=1, s>3/4$.
\vskip.1in
\begin{equation}
\begin{aligned}\label{k1}
\mu_{1,1}^T(u)=&\|u\|_{L_T^\infty H^s}+\|D_x^{s} \,\partial_x u\|_{L_x^\infty L_{yT}^2}+\|D_y^{s} \,\partial_x u\|_{L_x^\infty L_{yT}^2}\\
&+\|\partial_x u\|_{L_T^2 L_{xy}^\infty}+\|u\|_{L_x^{2}L_{yT}^{\infty}}.
\end{aligned}
\end{equation}

\item  For $k=2, s>3/4$.
\vskip.1in
\begin{equation}
\begin{aligned}\label{k2}
\mu_{1,2}^T(u)=&\|u\|_{L_T^\infty H^s}+\|D_x^{s} \,\partial_x u\|_{L_x^\infty L_{yT}^2}+\|D_y^{s} \,\partial_x u\|_{L_x^\infty L_{yT}^2}\\
&+\|u\|_{L_T^3 L_{xy}^\infty}+\|\partial_x u\|_{L_T^{\frac94} L_{xy}^\infty}+\| u\|_{L_x^{2}L_{yT}^{\infty}}.
\end{aligned}
\end{equation}

\item For $3\leq k\leq 7, s>3/4.$
\vskip.1in

\begin{equation}
\begin{aligned}\label{k7}
\mu_{1,k}^T(u)=&\|u\|_{L_T^\infty H^s}+\|D_x^{s} \,\partial_x u\|_{L_x^\infty L_{yT}^2}+\|D_y^{s} \,\partial_x u\|_{L_x^\infty L_{yT}^2}\\
&+\|u\|_{L_T^{p_k} L_{xy}^\infty}+\|\partial_x u\|_{L_T^{\frac{12}{5}} L_{xy}^\infty}+\| u\|_{L_x^{4}L_{yT}^{\infty}},
\end{aligned}
\end{equation}

where $p_k=\frac{12(k-1)}{7-12\gamma}$ and $\gamma\in (0,1/12).$
\vskip.1in
\item For $k\geq 8, s>s_k=1-2/k.$
\vskip.1in
\begin{equation}
\begin{aligned}\label{k8}
\mu_{1,k}^T(u)=&\|u\|_{L_T^\infty H^s}+\|D_x^{s} \,\partial_x u\|_{L_x^\infty L_{yT}^2}+\|D_y^{s} \,\partial_x u\|_{L_x^\infty L_{yT}^2}\\
&+\|\partial_x u\|_{L_x^\infty L_{yT}^2}+\|u\|_{L_T^{\frac{3k}{2}+} L_{xy}^\infty}+\|\partial_xu\|_{L_T^{\frac{3k}{k+2}} L_{xy}^\infty} +\|u\|_{L_x^{\frac{k}{2}}L_{yT}^{\infty}},
\end{aligned}
\end{equation}

\end{itemize}

The norms involved in these spaces reflect important properties of the associated group to the gZK equation, i.e. linear estimates  like the smoothing effect, Strichartz estimates, maximal function estimates, ... .
The standard argument is then carried out in the closed ball.

 $$B_a^T=\{u\in X_T;\;\; \mu_{1,k}^T(u)\leq a=2c\|u_0\|_{H^s(\R^2)}\},$$

of the metric space 

$$
X_T=\{u\in C([0,T]:H^s(\R^2));\;\;\mu_{1,k}^T(u)<\infty\}.
$$ 

By applying to the integral equation \eqref{in-inhoeqk} each norm in the definition of $\mu_{1,k}^T(u)$, linear estimates yield
\begin{equation}
\mu_{1,k}^T(u)\leq c\|u_0\|_{H^s}+cT^\gamma(\mu_{1,k}^T(u))^{k+1}
\end{equation}
where $\gamma$ is a positive constant.
From this point, the local existence time is chosen so that

\begin{equation}\label{ktimesize}
ca^kT^\gamma\leq 1/2,
\end{equation}

which implies that the time size  $T\sim \|u_0\|^{-\frac{k}{\gamma}}_{H^s}$ and that the local solution satisfies $\mu_{1,k}^T(u)\leq 2c\|u_0\|_{H^s}$.

\section{Proof of Theorem \ref{point-wise}}\label{S:3}

  We consider the unitary group of operators $\{W(t) : t\in\R\}$ in $L^2(\R^2)$ defined as
  \begin{equation}
  \label{def33}
  W(t)u_0(x,y)= ( e^{it (\xi^3 + \xi \eta^2)}\widehat{u}_0(\xi, \eta))^{\lor}(x,y).
  \end{equation}

  Thus, for $\alpha=r_1\in(0,1),$  \eqref{pro1} yields

  $$
  |x|^{r_1}\,W(t)u_0(x,y)= |x|^{r_1}( e^{it (\xi^3 + \xi \eta^2)}\widehat{u}_0(\xi, \eta))^{\lor}(x,y)= (D_{1,r_1} (e^{it (\xi^3 + \xi \eta^2)}\widehat{u}_0(\xi,\eta)))^{\lor}(x,y).
  $$
 and from Lemma \ref{Lemma1} that \begin{equation}
  \label{for1}
  D_{1,r_1} (e^{it(\xi^3 +  \xi \eta^2)}\,\widehat u_0)(\xi, \eta) = e^{it(\xi^3 + \xi \eta^2)}\, D_{1,r_1}\widehat u_0(\xi,\eta) + e^{it(\xi^3 + \xi \eta^2)} \Phi_{1, t,\alpha}(\widehat
  u_0)(\xi, \eta),
  \end{equation}
   with
  $$
  \| \Phi_{1,t,r_1}(\widehat u_0)\|_p\leq c_{r_1}(1+|t|)( \|\widehat u_0\|_p + \| \,(1+\xi^2 + \eta^2)^{r_1} \,\widehat u_0\|_p).
  $$
  Hence, taking Fourier transform in \eqref{for1}
  we obtain the identity
 \begin{equation}
  \label{for2}
 |x|^{r_1} \,W(t) u_0(x,y)= W(t)(|x|^{r_1} u_0)(x,y) + W(t) (\{ \Phi_{1,t,r_1}(\widehat u_0)(\xi,\eta)\}^{\lor})(x,y).
  \end{equation}
   with $\Phi_{1,t,r_1}$ as in \eqref{d10a} and
  \begin{equation}
  \label{estimate1}
 \begin{aligned}
  \| \{ \Phi_{1,t,r_1}(\widehat u_0)(\xi,\eta)\}^{\lor}\|_2&=\|\Phi_{1,t,r_1}(\widehat u_0)\|_2
 \\
 & \leq c_{r_1}(1+|t|) ( \|\widehat u_0\|_2 + \| \,(1+\xi^2 + \eta^2)^{r_1} \,\widehat u_0\|_2 \\
 & \leq c_{r_1}(1+|t|) (\|u_0\|_2 + \|\,D_x^{s}u_0\|_2 + \|\,D_y^{s}u_0\|_2).
 \end{aligned}
  \end{equation}

    On the other hand, if $\beta\in(0,r_1)$, then
\begin{equation}
  \label{for3}
	\begin{aligned}
D_x^{\beta}( |x|^{\alpha} \,W(t) u_0)(x,y)=& W(t)(D_x^{\beta}|x|^{r_1} u_0)(x,y) \\ 
&+ W(t) (D_x^{\beta}\{ \Phi_{1,t,r_1}(\widehat u_0)(\xi,\eta)\}^{\lor})(x,y).
  \end{aligned}
	\end{equation}

   In order to prove \eqref{d212-norm}  we need to show that
  \begin{equation}
  \label{est7}
  \begin{aligned}
 &  \| D_x^{\beta}
(\int \frac{e^{it(\varphi(\xi + \tau,\eta)- \varphi(\xi,\eta))}-1}
 {|\tau|^{1+r_1}}\,\widehat u_0(\xi+\tau, \eta)\, d\tau)^{\lor} \|_2
\\
\\
 &\leq c_{\alpha,\beta}(1+|t|) (\|u_0\|_2+\|D_x^{\beta+2\alpha}u_0\|_2 + \|D_y^{\beta+2\alpha}u_0\|_2).
  \end{aligned}
  \end{equation}
 Thus, we write
   \begin{equation}
  \label{est8}
  \begin{aligned}
 &  \| D^{\beta}_x
(\int \frac{e^{it(\varphi(\xi + \tau,\eta)- \varphi(\xi,\eta))}-1}
 {|\tau|^{1+r_1}}\,\widehat u_0(\xi+\tau, \eta)\, d\tau)^{\lor} \|_2 \\
& =\| \int \frac{|\xi|^{\beta}(e^{it(\varphi(\xi + \tau,\eta)- \varphi(\xi,\eta))}-1)}
 {|\tau|^{1+r_1}}\,\widehat u_0(\xi+\tau, \eta)\, d\tau \|_2\\
& \leq  \| \int \frac{|\xi|^{\beta}|e^{it(\varphi(\xi + \tau,\eta)- \varphi(\xi,\eta))}-1|}
 {|\tau|^{1+r_1}}\,|\widehat u_0(\xi+\tau, \eta)|\, d\tau \|_2\\
 & \leq c_{\beta} \| \int\frac{ |\xi+\tau|^{\beta}|e^{it(\varphi(\xi + \tau,\eta)- \varphi(\xi,\eta))}-1|}
 {|\tau|^{1+r_1}}\,|\widehat u_0(\xi+\tau, \eta)|\, d\tau \|_2\\
&\; \;\;+c_{\beta}  \| \int\frac{ |\tau|^{\beta}|e^{it(\varphi(\xi + \tau,\eta)- \varphi(\xi,\eta))}-1|}
 {|\tau|^{1+r_1}}\,|\widehat u_0(\xi+\tau, \eta)|\, d\tau \|_2\\
 &= I_1+I_2.
  \end{aligned}
 \end{equation}

   Following the argument used in the proof of Lemma \ref{Lemma1} to get \eqref{d10-norm} it holds that
   \begin{equation}
  \label{est9}
  \begin{aligned}
   I_1&\leq c_{r_1,\beta}(1+|t|)(\| |\xi|^{\beta} \widehat u_0\|_2
 +\|\,(1 + \xi^2 + \eta^2)^{r_1}|\xi|^{\beta} \widehat u_0\|_2)\\
   &\le c_{r_1,\beta}(1+|t|)(\| D_x^{\beta}u_0\|_2 + \| D_y^{\beta}u_0\|_2+\| D_x^{\beta+2r_1}u_0\|_2 + \| D_y^{\beta+2r_1}u_0\|_2)\\
	&\leq c_{r_1,\beta}(1+|t|)(\|u_0\|_2 + \|\,D_x^{\beta+s}u_0\|_2 + \|\,D_y^{\beta+s}u_0\|_2).
  \end{aligned}
 \end{equation}

 To bound $I_2$  we observe that this estimate is similar to that one used in the proof of Lemma \ref{Lemma1}
 with $r_1-\beta\geq 0$ instead of $\alpha$. Hence,
 \begin{equation}
  \label{est10}
  \begin{aligned}
  I_2 &\leq c_{r_1,\beta}(1+|t|)(\|  \widehat u_0\|_2 + \|\,(1 + \xi^2 + \eta^2)^{(r_1-\beta)} \widehat u_0\|_2)\\
  & \leq c_{r_1,\beta}(1+|t|)(\| u_0\|_2+\| D_x^{2(r_1-\beta)}u_0\|_2 + \| D_y^{2(r_1-\beta)}u_0\|_2)\\
	& \leq  c_{r_1,\beta}(1+|t|)(\|u_0\|_2 + \|\,D_x^{\beta+s}u_0\|_2 + \|\,D_y^{\beta+s}u_0\|_2).
  \end{aligned}
 \end{equation}

 For the weight in the $y$ direction the analysis follows similar arguments and hence we get that for $r_2\in(0,1)$  
 \begin{equation}
  |y|^{r_2} \,W(t) u_0(x,y)= W(t)(|y|^{r_2} u_0)(x,y) + W(t) (\{ \Phi_{2,t,r_2}(\widehat u_0)(\xi,\eta)\}^{\lor})(x,y).
  \end{equation}
   with $\Phi_{j,t,r_2}$ as in \eqref{d10b} and
  \begin{equation}
  \label{estimate3}
 \begin{aligned}
  \| \{ \Phi_{2,t,r_2}(\widehat u_0)(\xi,\eta)\}^{\lor}\|_2 \leq c_{r_2}(1+|t|) (\|u_0\|_2 + \|\,D_x^{s}u_0\|_2 + \|\,D_y^{s}u_0\|_2),
 \end{aligned}
  \end{equation}
    and similarly \eqref{d212-norm} for $j=2$ is obtained.
 
   This completes the proof of Theorem \ref{point-wise}.\\

\section{Proof of Theorem \ref{th3}}\label{S:4}

We consider the most interesting case $s=2\max\{r_1, r_2\},$ with $s_k<s<1$ as in the LWP theory in $H^s$. 

 \underline{Case 1:  $k=1.$}
\vskip.1in
From the previous assumption on  $s$, and with $ s>s_1=3/4.$ 

Let $u\in C([0,T] : H^{s}(\R))$ be  the unique solution of the ZK IVP 
satisfying the integral equation 

\begin{equation}
\label{in-inhoeq}
u(t)=\Psi(u(t))=W(t)u_0-\int_0^{t} W(t-t')(u \; \partial_x u)(t')dt',\;\;t\in[0,T]
\end{equation}

with $T=T(\|u_0\|_{H^s})<1$  in \eqref{ktimesize}  chosen to satisfy

\begin{equation}\label{timeprop}
c\mu_{1, 1}^T(u)T^\gamma\leq1/2,
\end{equation}
where $\gamma=1/2$ for $k=1$, $T\sim \|u_0\|_{H^s}^{-2}$  and $a=2c\|u_0\|_{H^s}$ is the radius of the ball in the contraction argument in $H^s$ so that 
\begin{equation}
\label{est-no}
\mu_{1,1}^T(u)\leq a=2c\|u_0\|_{H^s}.
\end{equation}
 
Now we suppose additionally that
 $$
 u_0\in Z_{s,(r_1,r_2)}=H^{}(\R^2)\cap L^2((\,|x|^{2r_1}+|y|^{2r_2})\,dxdy),
 $$
 and introduce the new norm
 \begin{equation}\label{wnorm}
\mu_{2, 1}^{T}(u) =\mu_{1,1}^{ T}(u) + \|\,(|x|^{r_1}+|y|^{r_2})u(t)\|_{L^{\infty}_{ T}L^2_{xy}}.
\end{equation}

Let us estimate $\mu_{2,1}^{T}(u) $ in the integral equation \eqref{in-inhoeq}  noticing that for the second term in the definition of $\mu_{2,1}^{T}(u) $ it is enough to consider
the norms

$$
 i)\;\;\;\|\,|x|^{r_1}u(t)\|_{L^{\infty}_{ T}L^2_{xy}} \;\;\; \text{and} \;\;\;ii)\;\;\; \|\,|y|^{r_2}u(t)\|_{L^{\infty}_{ T}L^2_{xy}}
$$

For $i)$ we have  from  Theorem \ref{point-wise} and \eqref{est-no} that

\begin{equation}
\label{est1c}
\begin{aligned}
\|\,|x|^{r_1}u\|_{L^2_{x,y}}&\leq \|\, |x|^{r_1} u_0 \|_2+ c(1+{T})\|u_0\|_{H^{2r_1}}+\|\,|x|^{r_1}\int_0^{t} W(t-t')(u \; \partial_x u)(t')dt'\|_{L^2_{xy}}\\
&\leq \|\, |x|^{r_1} u_0 \|_2+c(1+T)\|u_0\|_{H^{s}}+\int_0^{T}\;\||x|^{r_1}\;u\partial_xu\|_{L^2_{xy}}\;dt'\\
&+c(1+T)\int_0^{T}\|u\;\partial_xu\|_{H^{s}}\;dt'\\
&\leq \| \,|x|^{r_1} u_0 \|_2+c(1+T)\|u_0\|_{H^{s}}+c(1+T)T^\frac12(\mu_{1,1}^{T}(u))^2+I_{1},
\end{aligned}
\end{equation} 
where 
\begin{equation}
\label{kestim}
I_{1}=\int_0^{T}\;\|\,|x|^{r_1}\;u\partial_xu\|_{L^2_{xy}}\,dt'.
\end{equation}

 We obtain a similar estimate for $ii)$ and therefore we conclude that

\begin{equation}
\label{subtotalnorm}
\mu_{2,1}^{T}(u)\leq c(1+T)(\|u_0\|_{H^{s}}+\|\,(|x|^{r_1}+|y|^{r_2})u_0\|_{L^2_{x\,y}})+c(1+T)T^\frac12(\mu_{1,1}^{T}(u))^2+cI_{1}.
\end{equation}

Now we have for $I_1$:

\begin{equation}
\begin{aligned}
I_{1}&\leq \|\,|x|^{r_1}u\|_{L^{\infty}_{ T}L^2_{x\,y}}\|\partial_x u\|_{L^{1}_{ T}L^\infty_{xy}}\\
&\leq \mu_2^{T}(u)\;T^\frac12\|\partial_x u\|_{ L^2_{T}L^\infty_{xy}}\\
&\leq T^\frac12 \mu_{1,1}^{T}(u)\mu_{2,1}^{T}(u).
\end{aligned}
\end{equation}

In summary we get

\begin{equation}
\label{totalnorm}
\begin{aligned}
\mu_{2,1}^{T}(u)\leq &c(1+T)(\|u_0\|_{H^{s}}+\|\,(|x|^{r_1}+|y|^{r_2})u_0\|_{L^2_{xy}})\\
&+c(1+T)T^\frac12(\mu_{1,1}^{T}(u))^2+cT^\frac12\mu_{1,1}^{T}(u)\mu_{2,1}^{T}(u).
\end{aligned}
\end{equation}

From the time size in \eqref{timeprop}, we can pass the last term to the left side and obtain
\begin{equation}
\label{normbd}
\mu_{2,1}^{T}(u)\leq 2c(1+T)(\|u_0\|_{H^{s}}+\|\,(|x|^{r_1}+|y|^{r_2})u_0\|_{L^2_{x,y}})+2c(1+T)\|u_0\|_{H^{s}}.
\end{equation}

This basically completes the proof of this theorem in the $k=1$ case.
\vskip.1in
\underline{Case 2:  $k\geq2.$} 
\vskip.1in
For these nonlinearities the argument follows exactly the same ideas as in the former case and we provide some details at the points where the estimates depend on the norms involved in the associated local theory in $H^s$.

We again define a new norm as in \eqref{wnorm}
\begin{equation}\label{wnormk}
\mu_{2,k}^{T}(u) =\mu_{1,k}^{ T}(u) + \|\,(|x|^{r_1}+|y|^{r_2})u(t)\|_{L^{\infty}_{ T}L^2_{xy}},
\end{equation}
with $\mu^T_{1,k}$ the norm associated to the solution space $X_T$ in Theorem \ref{th1} and given  in section \ref{S:2}.

Now we consider   the second term in \eqref{wnormk} of the solution $u$ represented in Duhamel's formula and observe that it is enough consider the norms
$$
 i)\;\;\;\|\,|x|^{r_1}u(t)\|_{L^{\infty}_{ T}L^2_{xy}} \;\;\; \text{and} \;\;\;ii)\;\;\; \|\,|y|^{r_2}u(t)\|_{L^{\infty}_{ T}L^2_{xy}}.
$$
From Theorem \ref{point-wise} we can restrict our  attention to the first norm above and obtain:

\begin{equation}
\label{est1ck}
\begin{aligned}
\|\,|x|^{r_1}u\|_{L^2_{x,y}}&\leq \|\, |x|^{r_1} u_0 \|_2+ c(1+{T})\|u_0\|_{H^{2r_1}}+\|\,|x|^{r_1}\int_0^{t} W(t-t')(u^k \; \partial_x u)(t')dt'\|_{L^2_{xy}}\\
&\leq \|\, |x|^{r_1} u_0 \|_2+c(1+T)\|u_0\|_{H^{s}}+\int_0^{T}\;\|\,|x|^{r_1}\;u^k\partial_xu\|_{L^2_{xy}}\;dt'\\
&+c(1+T)\int_0^{T}\|u^k\;\partial_xu\|_{H^{s}}\;dt'\\
&\leq \|\, |x|^{r_1} u_0 \|_2+c(1+T)\|u_0\|_{H^{s}}+c(1+T)T^\gamma(\mu_{1,k}^{T}(u))^{k+1}+I_{k},
\end{aligned}
\end{equation} 
where 
\begin{equation}
\label{kestimk}
I_{k}=\int_0^{T}\;\|\,|x|^{r_1}\;u^k\partial_xu\|_{L^2_{xy}}\,dt'.
\end{equation}

Let us estimate $I_{k}$ for the different values of $k$:
\begin{itemize}
\item For $k=2$.

\begin{equation}
\begin{aligned}
I_{2}&\leq \|\,|x|^{r_1}u\|_{L^{\infty}_{ T}L^2_{x\,y}}\|u\partial_x u\|_{L^{1}_{ T}L^\infty_{xy}}\\
&\leq \mu_{2, 2}^{T}(u)\;\|u\|_{ L^\frac95_{T}L^\infty_{xy}}\|\partial_x u\|_{ L^\frac94_{T}L^\infty_{xy}}\\
&\leq T^\frac29\mu_{2, 2}^{T}(u)\|u\|_{ L^3_{T}L^\infty_{xy}}\|\partial_x u\|_{ L^\frac94_{T}L^\infty_{xy}}\\
&\leq T^\frac29(\mu_{1,2}^{T}(u))^2\mu_{2, 2}^{T}(u).
\end{aligned}
\end{equation}

\item For $3\leq k\leq 7.$

\begin{equation}
\begin{aligned}
I_{k}&\leq \|\,|x|^{r_1}u\|_{L^{\infty}_{ T}L^2_{xy}}\|u^{k-1}\partial_x u\|_{L^{1}_{ T}L^\infty_{xy}}\\
&\leq \mu_{2, k}^{T}(u)\|u^{k-1}\|_{ L^\frac{12}{7}_{T}L^\infty_{xy}}\|\partial_x u\|_{ L^\frac{12}{5}_{T}L^\infty_{xy}}\\
&\leq T^\gamma\mu_{2, k}^{T}(u)\|u\|^{k-1}_{ L^{\frac{12(k-1)}{7}}_{T}L^\infty_{xy}}\|\partial_x u\|_{ L^\frac{12}{5}_{T}L^\infty_{xy}}\\
&\leq T^\gamma(\mu_{1,k}^{T}(u))^{k}\mu_{2, k}^{T}(u).
\end{aligned}
\end{equation}

\item For $k\geq 8.$

\begin{equation}
\begin{aligned}
I_{k}&\leq \|\,|x|^{r_1}u\|_{L^{\infty}_{ T}L^2_{xy}}\|u^{k-1}\partial_x u\|_{L^{1}_{ T}L^\infty_{xy}}\\
&\leq \mu_{2, k}^{T}(u)\|u^{k-1}\|_{ L^\frac{3k}{2(k-1)}_{T}L^\infty_{xy}}\|\partial_x u\|_{ L^\frac{3k}{k+2}_{T}L^\infty_{xy}}\\
&\leq \mu_{2, k}^{T}(u)\|u\|^{k-1}_{ L^{\frac{3k}{2}}_{T}L^\infty_{xy}}\|\partial_x u\|_{ L^\frac{3k}{k+2}_{T}L^\infty_{xy}}\\
&\leq T^\gamma\mu_{2, k}^{T}(u)\|u\|^{k-1}_{ L^{\frac{3k}{2} +}_{T}L^\infty_{xy}}\|\partial_x u\|_{ L^\frac{3k}{k+2}_{T}L^\infty_{xy}}\\
&\leq T^\gamma(\mu_{1,k}^{T}(u))^{k}\mu_{2, k}^{T}(u).
\end{aligned}
\end{equation}

\end{itemize}

Hence the following estimate for the norm in \eqref{wnormk} holds

\begin{equation}
\label{totalnormk}
\begin{aligned}
\mu_{2,k}^{T}(u)\leq &c(1+T)(\|u_0\|_{H^{s}}+\|\,(|x|^{r_1}+|y|^{r_2})u_0\|_{L^2_{xy}})\\
&+c(1+T)T^\gamma(\mu_{1,k}^{T}(u))^{k+1}+cT^\gamma(\mu_{1,k}^{T}(u))^{k}\mu_{2, k} ^{T}(u).
\end{aligned}
\end{equation}

From the time size in in the contraction argument in Theorem \ref{th1}, we can again pass the last term to the left side and obtain
\begin{equation}
\label{normbdk}
\mu_{2,k}^{T}(u)\leq 2c(1+T)(\|u_0\|_{H^{s}}+\|\,(|x|^{r_1}+|y|^{r_2})\,u_0\|_{L^2_{xy}})+2c(1+T) \|u_0\|_{H^{s}}.
\end{equation}

Which completes the proof of the theorem.




\end{document}